\documentclass[11pt,a4paper]{article}
\usepackage[T1]{fontenc}
\usepackage{lmodern}
\usepackage{amsmath,amssymb,amsthm}
\usepackage[left=25mm,right=25mm,top=24mm,bottom=25mm]{geometry}
\usepackage{enumitem}
\usepackage{tikz}
\usepackage[hidelinks]{hyperref}
\numberwithin{equation}{section}
\newtheorem{theorem}{Theorem}[section]
\newtheorem{lemma}[theorem]{Lemma}
\newtheorem{proposition}[theorem]{Proposition}
\newtheorem{corollary}[theorem]{Corollary}
\theoremstyle{definition}

\newtheorem{conjecture}[theorem]{Conjecture}

\DeclareMathOperator{\Def}{Def}
\DeclareMathOperator{\PR}{PR}

\newcommand{\IS}{\mathcal I}
\newcommand{\cf}{\chi_f}
\hypersetup{pdftitle={Independent Sets in C4-Free Planar Graphs: Sharp Claw-Free Bounds and Structural Reductions},pdfsubject={Weighted independence bounds, residual matchings, and the five-sixteenths conjecture}}
\title{Independent Sets in $C_4$-Free Planar Graphs:\\ Sharp Claw-Free Bounds and Structural Reductions}
\author{Zhou Guofei\footnote{Email: gfzhou@nju.edu.cn}
\\
{\small School of Mathematics, Nanjing University, Nanjing 210093, China}
 }
\date{}
\begin{document}
\maketitle
\vspace{-1.7em}
\begin{abstract}
For every claw-free $C_4$-free planar graph, we prove the sharp weighted bound $\alpha_w(G)\ge (10/31)\sum_v w(v)$, equivalently $\chi_f(G)\le31/10$. The proof uses a line-graph representation and the matching polytope; the exceptional small roots are excluded analytically by planar duality and degree-four neighborhood analysis in triangulations. We also prove that, in every $2$-connected $C_4$-free graph of minimum degree at least four, an independent selection of claw centers leaves a compatible subgraph of noncenters with fractional chromatic number at most three; planarity is not required. We apply these results to the conjecture $\alpha(G)\ge5|V(G)|/16$. A minimum-order counterexample without degree-three vertices has at least forty-nine vertices. When all degree-three components of a minimum counterexample are triangle trees, we derive exact boundary states and a joint matching formula. We identify the remaining deficiency estimate and the additional component types that must be addressed to extend this approach. The general conjecture remains open.
\end{abstract}


\section{Introduction}
All graphs in this paper are finite, simple, and undirected. A graph is $C_4$-free if it has no four-cycle as a subgraph; the exclusion is not restricted to induced cycles. The independence number $\alpha(G)$ is the maximum cardinality of an independent set. We study the following conjecture.
\begin{conjecture}\label{conj:main}
Every $n$-vertex $C_4$-free planar graph satisfies
\begin{equation}\label{eq:target}
\alpha(G)\ge\frac{5n}{16}.
\end{equation}
\end{conjecture}
Kelly, Kolichala, McFarland, and Su~\cite{Kelly2026} prove $\alpha(G)\ge30n/119$ for every $C_4$-free planar graph and ask for the supremum of $|V(G)|/\alpha(G)$ in this class. Their result follows from a theorem for the larger class of planar graphs with no triangle sharing an edge with a four-cycle. Their Figure~1 gives a sixteen-vertex $C_4$-free planar graph with independence number five, so the coefficient in~\eqref{eq:target}, if valid, is best possible. The conjecture also has an earlier planar Ramsey formulation due to Sun, Yang, Lin, and Song~\cite{Sun2008}; we establish the equivalence in Proposition~\ref{prop:ramsey}.

\subsection{Main results and approach}
Our first result determines the optimal weighted independence coefficient for the claw-free subclass. For nonnegative vertex weights $w$, write $w(X)=\sum_{v\in X}w(v)$ and let $\alpha_w(G)$ be the maximum weight of an independent set. We write $\chi_f(G)$ for the fractional chromatic number.
\begin{theorem}\label{thm:main-claw}
Every claw-free $C_4$-free planar graph $H$, with nonnegative vertex weights $w$, satisfies
\[
\cf(H)\le\frac{31}{10},\qquad
\alpha_w(H)\ge\frac{10}{31}\sum_{v\in V(H)}w(v).
\]
In particular, $\alpha(H)\ge10|V(H)|/31$. Both constants are sharp.
\end{theorem}
Theorem~\ref{thm:main-claw} concerns a narrower graph class than the bound of Kelly et al., but gives a sharp conclusion for every nonnegative weighting, rather than only for the independence number. Its proof uses a line-graph representation: independent sets become matchings in a planar subcubic root of girth at least five. Equivalently, the root has fractional edge-chromatic number at most $31/10$. Edmonds' matching-polytope theorem~\cite{Edmonds1965} reduces the bound to odd-set density. Planar duality reduces the two smallest possible exceptional root orders to three triangulation configurations, which we exclude by local degree arguments and the nonplanarity of $K_{3,3}$. Section~\ref{sec:sharp} gives this analytic proof and a sharp example. No computational enumeration is used in this proof.

A closely related result is due to Gimbel, K\"undgen, Li, and Thomassen~\cite[Lemma~3.4]{Gimbel2019}. For planar $2$-edge-connected graphs of maximum degree three, they identify the unique-degree-two-vertex case as the exception to a fractional edge-chromatic bound of three. Their proof uses the same classical odd-set-density characterization: if an exceptional root has order $u$, its density is $(3u-1)/(u-1)=3+2/(u-1)$. The contribution of Theorem~\ref{thm:root-sharp} is to impose girth at least five, exclude the two smallest possible exceptional orders, and obtain the sharp bound $31/10$ for all planar subcubic roots of that girth, including roots with bridges or several components. We do not claim the matching-polytope method or the unique-degree-two obstruction as new.

Our second principal result allows an independent set to be chosen among claw centers before applying a matching bound to the remaining vertices. A \emph{claw} is an induced $K_{1,3}$, and its center is the vertex of degree three within that induced subgraph.
\begin{theorem}\label{thm:main-residual}
Let $G$ be a $2$-connected $C_4$-free graph with $\delta(G)\ge4$. Let $B$ be its claw centers and $A=V(G)\setminus B$. For every independent set $I\subseteq B$ and every nonnegative vertex weighting,
\begin{equation}\label{eq:main-residual}
\alpha_w(G)\ge w(I)+\frac13w(A\setminus N_G(I)).
\end{equation}
Indeed, $\cf(G[A\setminus N_G(I)])\le3$.
\end{theorem}
Unlike Theorem~\ref{thm:main-claw}, this statement does not require planarity. The proof identifies the odd-set obstructions to the uniform matching-polytope vector $x_e=1/3$ in a residual triangle root. Each possible obstruction would produce a cutvertex in the original graph. The residual root itself need not be $2$-edge-connected, so the bridgeless-root result of~\cite[Lemma~3.4]{Gimbel2019} cannot simply be applied to it. The additional argument controls root obstructions through their attachments in the original graph; the absence of a planarity assumption alone is not the distinction, since the odd-set argument in that lemma does not use planarity either.

We then apply the two matching bounds to a hypothetical counterexample of minimum order to Conjecture~\ref{conj:main}. Deletion stability and the absence of clique cutsets yield restrictions on its degree-three vertices. In particular, a minimum counterexample with no degree-three vertex has order at least forty-nine (Theorem~\ref{thm:no-three-order}). If every degree-three component consists of disjoint triangles joined in a tree by bridges, we obtain exact boundary states and a joint optimization formula coupling those components with the residual root matchings (Theorem~\ref{thm:joint}). These are necessary conditions on a minimum counterexample, not a proof of the general conjecture.

\paragraph{Organization.}
Section~\ref{sec:prelim} records the local structure and matching-polytope facts. Sections~\ref{sec:sharp} and~\ref{sec:residual} prove the two main theorems. Section~\ref{sec:minimal} develops the minimum-counterexample framework, and Section~\ref{sec:applications} gives its minimum-degree-four applications. Sections~\ref{sec:branches} and~\ref{sec:joint} treat triangle-tree boundaries and the joint matching model. Section~\ref{sec:next} states the remaining obstacles; Appendix~\ref{app:illustrations} illustrates why adjacency of the degree-four vertices matters in the small-root exclusion.

\subsection{Notation and the Ramsey reformulation}
For a vertex set $X$, let $G-X=G[V(G)\setminus X]$, let $N_G(X)$ denote its open neighborhood, and let $N_G[X]=X\cup N_G(X)$. We use $N_G(X)\setminus X$ when an external boundary is intended. The family of independent sets of $G$ is $\IS(G)$. A plane graph has a fixed embedding, and face length counts edge occurrences in its boundary walk. Write $\nu(R)$ for the matching number and
\[
\Def(R)=|V(R)|-2\nu(R).
\]
For nonnegative weights, $\alpha_w$ and $\nu_w$ denote the corresponding maximum total weights. The fractional chromatic number assigns nonnegative weights to independent sets, covering each vertex with total weight at least one, while minimizing total weight.

The planar Ramsey number $\PR(C_4,K_t)$ is the least $N$ such that every $C_4$-free planar graph of order $N$ has an independent set of size $t$. Sun, Yang, Lin, and Song~\cite{Sun2008} conjectured the formula in the next proposition; Chen, Cheng, Zhang, and Zhou~\cite{Chen2014} established its case $t=8$ and recorded the earlier cases $t\le7$.
\begin{proposition}\label{prop:ramsey}
Conjecture~\ref{conj:main} is equivalent to the identities
\begin{equation}\label{eq:ramsey}
\PR(C_4,K_t)=\left\lfloor\frac{16(t-1)}5\right\rfloor+1
=3t+\left\lfloor\frac{t-1}{5}\right\rfloor-2
\end{equation}
for every integer $t\ge2$.
\end{proposition}
\begin{proof}
Let $P$ be a sixteen-vertex $C_4$-free planar graph with $\alpha(P)=5$. Write $t-1=5q+r$, where $0\le r\le4$. The disjoint union of $q$ copies of $P$ and $r$ triangles has order $16q+3r=\lfloor16(t-1)/5\rfloor$ and independence number $t-1$. This proves the lower bound in~\eqref{eq:ramsey}. Under~\eqref{eq:target}, a graph of the next order has independence number greater than $t-1$, proving the upper bound.

Conversely, assume~\eqref{eq:ramsey} for all $t$ and put $a=\alpha(G)$. If $|V(G)|>\lfloor16a/5\rfloor$, an induced subgraph of order $\lfloor16a/5\rfloor+1$ would have an independent set of size $a+1$, a contradiction. The empty graph is immediate.
\end{proof}
The established cases $t\le8$ imply that any counterexample has independence number at least eight and order at least twenty-six.

\section{Neighborhood structure and fractional matchings}\label{sec:prelim}
\begin{lemma}\label{lem:neighbors}
In a $C_4$-free graph, two distinct vertices have at most one common neighbor. Every open neighborhood induces a matching and isolated vertices. Every clique has size at most three, and every edge belongs to exactly one maximal clique.
\end{lemma}
\begin{proof}
Two common neighbors give a four-cycle. A two-edge path in $N(v)$ gives a four-cycle with $v$, so the neighborhood has maximum degree at most one. A four-clique contains a four-cycle. Two distinct triangles containing the same edge also give a four-cycle, proving uniqueness of the maximal clique containing each edge.
\end{proof}

We use Edmonds' description of the matching polytope~\cite{Edmonds1965}: it is the convex hull of matching incidence vectors, characterized by
\begin{align}
x_e&\ge0,\label{eq:mp1}\\
\sum_{e\ni v}x_e&\le1\quad(v\in V(R)),\label{eq:mp2}\\
\sum_{e\in E(R[U])}x_e&\le\frac{|U|-1}{2}\quad(|U|\text{ odd}).\label{eq:mp3}
\end{align}
If a uniform vector $x_e=c>0$ belongs to this polytope, then
\begin{equation}\label{eq:uniform}
\nu_w(R)\ge c\sum_{e\in E(R)}w(e),\qquad \cf(L(R))\le1/c.
\end{equation}
Indeed, the vector is a convex combination of matchings. Its weighted average is $c\sum_e w(e)$, and scaling that combination by $1/c$ gives a fractional coloring of the line graph.

For a subcubic graph $R$, define
\[
P_R(U)=\sum_{v\in U}(3-d_R(v)),\qquad r_R(U)=|\delta_R(U)|.
\]
We retain isolated vertices when forming residual roots. The next lemma records the standard odd-set argument in a form that also accounts for bridges and disconnected roots; compare~\cite[Lemma~3.4]{Gimbel2019}.
\begin{lemma}\label{lem:third}
The uniform vector $x_e=1/3$ belongs to the matching polytope of a subcubic graph $R$ if and only if
\begin{equation}\label{eq:ports}
P_R(U)+r_R(U)\ge3\qquad(|U|\text{ odd}).
\end{equation}
Equivalently, neither of the following occurs: an odd component with total degree deficit one; a bridge with an odd side all of whose vertices have degree three in $R$. If the vector is feasible, then
\begin{equation}\label{eq:def-ports}
\nu(R)\ge\left\lceil\frac{|E(R)|}{3}\right\rceil,
\qquad \Def(R)\le\frac{P_R(V(R))}{3}.
\end{equation}
\end{lemma}
\begin{proof}
The identity
\[
2e_R(U)=3|U|-P_R(U)-r_R(U)
\]
makes~\eqref{eq:ports} equivalent to~\eqref{eq:mp3}. The vertex constraints hold automatically. For odd $U$, the sum $P_R(U)+r_R(U)$ is positive and odd, so failure means that it is one. Splitting $R[U]$ into its components, one odd component has budget one, while any component of budget zero is cubic and has even order. The component of budget one has either deficit one and no exit, or deficit zero and one exit. These are exactly the two stated obstructions. The final bounds follow from~\eqref{eq:uniform} and $2|E(R)|=3|V(R)|-P_R(V(R))$.
\end{proof}

We will also use the following elementary observation about plane embeddings.
It does not require the embedded graph to be connected.
\begin{lemma}[A connected witness for cofacial ports]\label{lem:cofacial-witness}
Let $H$ be a plane graph and let $J$ be a connected graph drawn disjointly
from $H$ in the same plane. Suppose that for each $p$ in a set
$P\subseteq V(H)$ there is an arc from $p$ to a point of $J$ whose only
intersection with $H$ is $p$. Then all vertices of $P$ are incident with
one face of $H$. For each component $H'$ of $H$, the vertices of
$P\cap V(H')$ are likewise incident with one face of the inherited
embedding of $H'$.
\end{lemma}
\begin{proof}
The connected drawing of $J$ is contained in one component $f$ of the
complement of the drawing of $H$, that is, in one face of $H$. Each of
the specified arcs, with its endpoint $p$ removed, is a connected subset
of that complement meeting $J$, so it lies in $f$. Thus $p$ is incident
with $f$. The face $f$ is contained in a face of every component $H'$
considered on its own, proving the final assertion.
\end{proof}

\section{A sharp bound for the claw-free subclass}\label{sec:sharp}
\begin{lemma}\label{lem:line-root}
Every claw-free $C_4$-free planar graph $H$ is the line graph of a simple planar subcubic graph $R$ of girth at least five. Forests are allowed.
\end{lemma}
\begin{proof}
Use the maximal cliques of $H$, including singleton cliques for isolated vertices. By Lemma~\ref{lem:neighbors} and claw-freeness, each vertex belongs to at most two maximal cliques. Make a root vertex for each clique. Represent a vertex belonging to two cliques by an edge between their root vertices; for one clique use a new private leaf. Two cliques cannot share two vertices, so the root is simple. A clique vertex has degree equal to its clique size, at most three, and private leaves have degree one. Adjacency in $H$ is exactly intersection of the representing edges, hence $H=L(R)$.

A triangle in $R$ would represent three vertices forming a triangle in $H$. Their unique maximal triangle clique would force their three root edges to have a common endpoint, impossible for a root triangle. A root four-cycle gives a four-cycle in $H$. Finally, the planar line-graph criterion of Sedl\'a\v cek~\cite[Theorem~1]{Sedlacek1990} implies that $R$ is planar.
\end{proof}

\begin{lemma}[Small-root exclusion]\label{lem:finite}
There is no simple planar graph of girth at least five on seventeen or nineteen vertices with exactly one vertex of degree two and every other vertex of degree three.
\end{lemma}
\begin{proof}
Let $R$ be such a graph and let $x$ be its degree-two vertex. The planar girth bound gives at least twenty vertices for any cubic component. Thus $R$ is connected.

For a connected cyclic subgraph on $u$ vertices with total internal degree deficit $p$, the degree identity and the planar girth bound give
\begin{equation}\label{eq:side-size}
\frac{3u-p}{2}\le\frac{5(u-2)}3,
\qquad u\ge20-3p.
\end{equation}
Suppose an edge is a bridge. A tree side with $p_i$ original degree-two vertices would satisfy $3u-p_i-1=2(u-1)$, or $u=p_i-1$, impossible because $p_i\le1$. Both sides are cyclic. Their internal deficits are $p_i+1$, so~\eqref{eq:side-size} gives total order at least $34-3=31$. Hence $R$ has no bridge. A bridgeless connected subcubic graph has no cutvertex, so $R$ is $2$-connected.

For a two-edge bond, two cyclic sides similarly force total order at least $28-3=25$. A tree side instead satisfies $3u-p_i-2=2(u-1)$, or $u=p_i$. Therefore every two-edge bond isolates $x$.

The plane dual has no loops and exactly one parallel pair, corresponding to the two edges incident with $x$. Deleting one edge of this pair gives a simple triangulation $T$. Its two endpoints remain adjacent. For order seventeen, $R$ has twenty-five edges and ten faces; every face is pentagonal. Thus $T$ has degree sequence $4^2 5^8$, with the two degree-four vertices adjacent. For order nineteen, $R$ has twenty-eight edges and eleven faces, ten pentagons and one hexagon. The degree sequence of $T$ is either $4^1 5^{10}$ or $4^2 5^8 6^1$; in the second case the two degree-four vertices must be adjacent.

We exclude these three triangulations by examining the neighborhood of a
vertex of degree four. All embeddings below are regarded as embeddings on
the sphere, so either side of a cycle may be chosen as its interior.

\emph{Triangular regions.}
A side of a triangular cycle that contains at least one vertex, all of whose
interior vertices have degree at least five in $T$, contains at least four
vertices. Indeed, an interior vertex can have neighbors only among the
other interior vertices and the three boundary vertices. With one or two
interior vertices its degree is at most four. With three interior vertices,
each must be adjacent to all three boundary vertices in order to have
degree at least five. These edges contain a $K_{3,3}$, contradicting
planarity.

\emph{The neighborhood of a degree-four vertex.}
Choose a degree-four vertex $v$; when there are two such vertices, let the
other one be a neighbor $x_0$ of $v$. List the neighbors of $v$ in cyclic
order as $x_0,x_1,x_2,x_3$. The four faces at $v$ form a disk bounded by
$C=x_0x_1x_2x_3x_0$. Set
\[
W=V(T)\setminus\bigl(\{v\}\cup V(C)\bigr).
\]
Thus $|W|$ is five or six, and every vertex of $W$ has degree at least five.
The cycle $C$ has no chord. For if $x_0x_2$ were a chord, it would lie on
the side of $C$ not containing $v$ and split that disk into two triangular
regions. Both regions must contain a vertex: otherwise $x_1$ or $x_3$
would have only its two neighbors on $C$ and $v$, contrary to minimum degree
four. The triangular-region observation would then give $|W|\ge8$.
The other possible chord is excluded in the same way.

For each $i$ modulo four, let $z_i\in W$ be the third vertex of the face
incident with $x_ix_{i+1}$ on the side not containing $v$. In the rotation
at $x_i$, the neighbors on that side form a sequence of distinct vertices
of $W$, beginning at $z_{i-1}$ and ending at $z_i$. Its length is
$d_T(x_i)-3$, and consecutive terms are adjacent. In particular,
\begin{itemize}[nosep]
\item if $d_T(x_i)=4$, then $z_{i-1}=z_i$;
\item if $d_T(x_i)=5$, then $z_{i-1}\ne z_i$ and $z_{i-1}z_i\in E(T)$;
\item if $d_T(x_i)=6$, the sequence is $z_{i-1},t,z_i$ for three
      distinct vertices, with edges $z_{i-1}t$ and $tz_i$.
\end{itemize}
Once these sequences are specified, the neighborhoods of all four $x_i$
are exhausted. Every vertex of $W$ not appearing in a sequence therefore
has all its neighbors in $W$. Figure~\ref{fig:four-neighborhood} shows the
local configuration, including the identification used in Case~2.

\begin{figure}[tbp]
\centering
\begin{minipage}[t]{0.48\textwidth}
\centering
\begin{tikzpicture}[scale=.83,line cap=round,line join=round,
  every node/.style={font=\small},dot/.style={circle,fill=black,inner sep=1.5pt}]
\coordinate (v) at (0,0);
\coordinate (x0) at (-1,1);
\coordinate (x1) at (1,1);
\coordinate (x2) at (1,-1);
\coordinate (x3) at (-1,-1);
\coordinate (z0) at (0,2.5);
\coordinate (z1) at (2.5,0);
\coordinate (z2) at (0,-2.5);
\coordinate (z3) at (-2.5,0);
\foreach \i in {0,...,3} {\draw[line width=.45pt] (v)--(x\i);}
\draw[line width=.9pt] (x0)--(x1)--(x2)--(x3)--cycle;
\draw[line width=.45pt] (x0)--(z0)--(x1)--(z1)--(x2)--(z2)--(x3)--(z3)--(x0);
\node[dot,label=above:$v$] at (v) {};
\node[dot,label=left:$x_0$] at (x0) {};
\node[dot,label=right:$x_1$] at (x1) {};
\node[dot,label=right:$x_2$] at (x2) {};
\node[dot,label=left:$x_3$] at (x3) {};
\node[dot,label=above:$z_0$] at (z0) {};
\node[dot,label=right:$z_1$] at (z1) {};
\node[dot,label=below:$z_2$] at (z2) {};
\node[dot,label=left:$z_3$] at (z3) {};
\end{tikzpicture}

\small (a) Four distinct exterior face vertices.
\end{minipage}\hfill
\begin{minipage}[t]{0.48\textwidth}
\centering
\begin{tikzpicture}[scale=.83,line cap=round,line join=round,
  every node/.style={font=\small},dot/.style={circle,fill=black,inner sep=1.5pt}]
\coordinate (v) at (0,0);
\coordinate (x0) at (-1,1);
\coordinate (x1) at (1,1);
\coordinate (x2) at (1,-1);
\coordinate (x3) at (-1,-1);
\coordinate (a) at (-2.5,0);
\coordinate (b) at (2.5,0);
\coordinate (c) at (0,-2.5);
\foreach \i in {0,...,3} {\draw[line width=.45pt] (v)--(x\i);}
\draw[line width=.9pt] (x0)--(x1)--(x2)--(x3)--cycle;
\draw[line width=.45pt] (x0)--(a)--(x3)--(c)--(x2)--(b)--(x1);
\draw[line width=.45pt] (a).. controls (-2.2,2.65) and (.2,2.65)..(x1);
\node[dot,label=above:$v$] at (v) {};
\node[dot,label=above:$x_0$] at (x0) {};
\node[dot,label=right:$x_1$] at (x1) {};
\node[dot,label=right:$x_2$] at (x2) {};
\node[dot,label=left:$x_3$] at (x3) {};
\node[dot,label=left:$a$] at (a) {};
\node[dot,label=right:$b$] at (b) {};
\node[dot,label=below:$c$] at (c) {};
\path[use as bounding box] (-3.1,-3.0) rectangle (3.1,3.0);
\end{tikzpicture}

\small (b) $z_3=z_0=a$, $z_1=b$, $z_2=c$.
\end{minipage}
\caption{The degree-four neighborhood in Lemma~\ref{lem:finite}.
The thick edges form $C=x_0x_1x_2x_3x_0$; the four triangular faces at $v$
and the face across each edge of $C$ are shown. Panel (b) has
$d_T(x_0)=4$. These are partial local drawings, not complete candidate
triangulations: other exterior edges and vertices, including a possible
middle neighbor at a degree-six vertex, are omitted.}
\label{fig:four-neighborhood}
\end{figure}
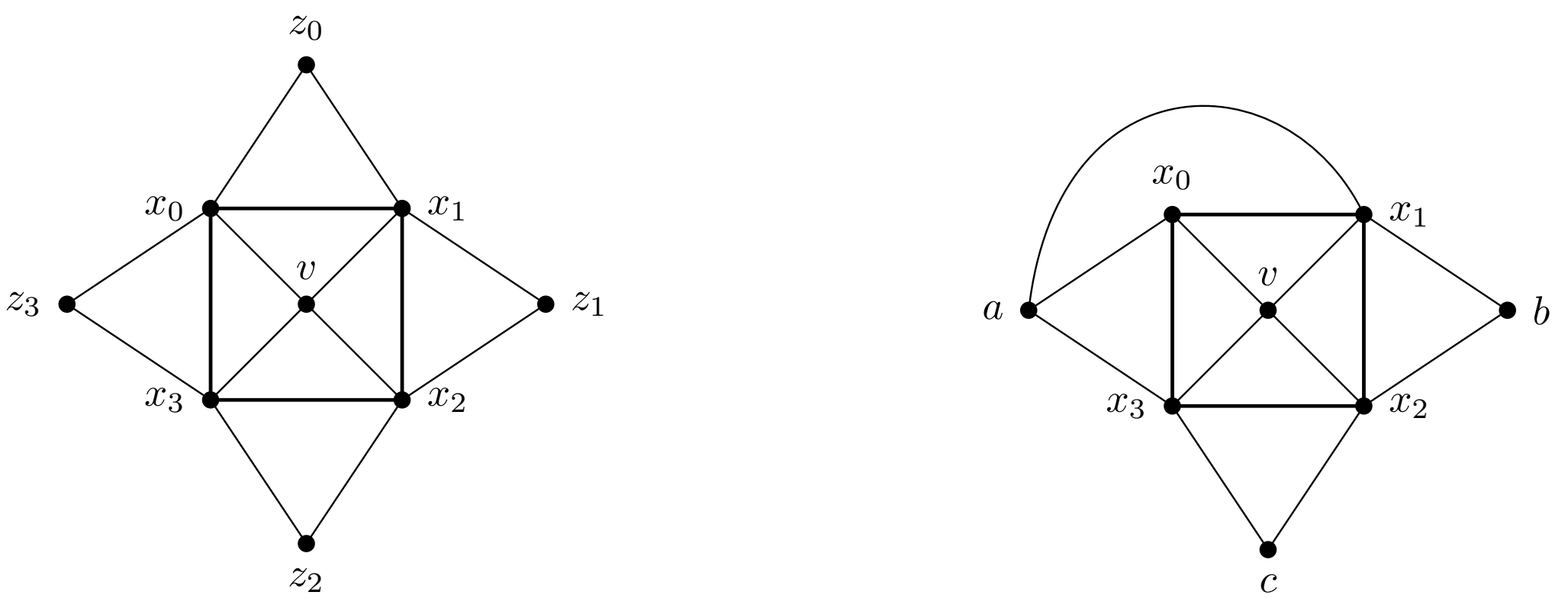

\emph{Case 1: degree sequence $4^1 5^{10}$.}
All four $x_i$ have degree five. We claim that the four $z_i$ are distinct.
Consecutive ones are distinct by the rotation rule. If $z_0=z_2$, this
vertex is adjacent to all four $x_i$ and to $z_1,z_3$. When $z_1\ne z_3$
this gives degree at least six. When $z_1=z_3$, the three distinct vertices
$v,z_0,z_1$ are each adjacent to all four $x_i$, again giving a $K_{3,3}$.
The identification $z_1=z_3$ is excluded symmetrically.

The two vertices $p,q$ of $W\setminus\{z_0,z_1,z_2,z_3\}$ have degree
five and all their neighbors in the six-vertex set $W$. Each is therefore
adjacent to every other vertex of $W$. But every $z_i$ then has the six
distinct neighbors
$x_i,x_{i+1},z_{i-1},z_{i+1},p,q$, a contradiction.

\emph{Case 2: two adjacent degree-four vertices.}
Here $x_0$ has degree four, so $z_3=z_0=:a$. Put $b=z_1$ and $c=z_2$.
The rotation rules at $x_1,x_2,x_3$, whose degrees are at least five, imply
that $a,b,c$ are pairwise distinct.

If all three of $x_1,x_2,x_3$ have degree five, then $ab,bc,ca$ are edges,
and $a$ already has the five distinct neighbors $x_0,x_1,x_3,b,c$.
Every vertex in $W\setminus\{a,b,c\}$ has all its neighbors in $W$.
For $|W|=5$ this is incompatible with degree at least five.
For $|W|=6$, the three remaining vertices must each be adjacent to all
other vertices of $W$, giving $d_T(a)\ge8>\Delta(T)$.

It remains to consider $|W|=6$ with one of $x_1,x_2,x_3$ of degree six.
All vertices of $W$ now have degree five. Let $t$ be the middle vertex
in the exterior neighbor sequence of the degree-six vertex.
If $t\notin\{a,b,c\}$, the two remaining vertices of $W$ are each adjacent
to every other vertex of $W$. Independently of the position of the
degree-six vertex, $a$ is adjacent to $x_0,x_1,x_3$ and at least two of
$b,c,t$. The two remaining vertices give $d_T(a)\ge7$, a contradiction.

If $t\in\{a,b,c\}$, distinctness within the exterior sequence leaves
only three possibilities. When $x_1$ has degree six, $t=c$, and $c$ already
has the five distinct neighbors $x_1,x_2,x_3,a,b$. The three vertices in
$W\setminus\{a,b,c\}$ must each be adjacent to every other vertex of $W$,
giving $d_T(c)\ge8$. The case where $x_3$ has degree six is symmetric,
with $t=b$. Finally, if $x_2$ has degree six, then $t=a$ and $a$ has the
six distinct neighbors $x_0,x_1,x_2,x_3,b,c$, contradicting its degree five.
This excludes all three triangulations and completes the proof.
\end{proof}

\begin{theorem}\label{thm:root-sharp}
Every simple planar subcubic graph $R$ of girth at least five supports the uniform fractional matching $x_e=10/31$. Hence
\[
\nu_w(R)\ge\frac{10}{31}\sum_{e\in E(R)}w(e),
\qquad\cf(L(R))\le\frac{31}{10}.
\]
\end{theorem}
\begin{proof}
The vertex constraints hold since $30/31<1$. For an odd set $U$ of size $u$, put $p=3u-2e_R(U)$. This is a positive odd integer. If $p\ge3$, then $e_R(U)\le3(u-1)/2$, which implies~\eqref{eq:mp3} because $10/31<1/3$.

If $p=1$, the induced graph has one degree-two vertex and all other vertices of degree three. Its degree-two component has at least seventeen vertices by~\eqref{eq:side-size}, and Lemma~\ref{lem:finite} excludes total orders seventeen and nineteen. Consequently $u\ge21$, and
\[
\frac{10}{31}e_R(U)=\frac{10}{31}\frac{3u-1}{2}
\le\frac{u-1}{2}.
\]
Thus all matching-polytope constraints hold, and~\eqref{eq:uniform} applies.
\end{proof}

\begin{proof}[Proof of Theorem~\ref{thm:main-claw} and sharpness]
Apply Lemma~\ref{lem:line-root} and Theorem~\ref{thm:root-sharp}. To see sharpness, subdivide one edge of the dodecahedral graph once, obtaining $R_0$ with twenty-one vertices and thirty-one edges. It remains planar, subcubic, and of girth five. Odd order gives $\nu(R_0)\le10$, and Theorem~\ref{thm:root-sharp} gives the reverse inequality. The line graph $H_0=L(R_0)$ is planar and claw-free. It is $C_4$-free: four mutually cyclically adjacent root edges would force a root four-cycle, a triangle with an additional incident edge, or a vertex of degree at least four; none is possible. Therefore $|V(H_0)|=31$ and $\alpha(H_0)=10$. The inequality $|V(H_0)|/\alpha(H_0)\le\cf(H_0)$ proves equality in the fractional bound. Disjoint copies prove optimality of the independence coefficient.
\end{proof}

\begin{corollary}\label{cor:claw-delete}
If $G$ is $C_4$-free and planar and $G-S$ is claw-free, then
\[
\alpha(G)\ge\left\lceil\frac{10(|V(G)|-|S|)}{31}\right\rceil.
\]
In particular, $|S|\le|V(G)|/32$ implies~\eqref{eq:target}.
\end{corollary}

\section{Residual triangle roots}\label{sec:residual}
In this section $G$ is any $2$-connected $C_4$-free graph with minimum degree at least four. Let $B$ be its claw centers and $A=V(G)\setminus B$. By Lemma~\ref{lem:neighbors}, every $a\in A$ has degree four and belongs to exactly two triangles, which contain all its incident edges.

Let $\mathcal T$ be the triangles of $G$. Form a graph $R$ with vertex set $\mathcal T$ and an edge for each $a\in A$, joining its two triangles. The root is simple and subcubic. For $I\subseteq B$, let $C(I)\subseteq\mathcal T$ be the triangles meeting $I$, and put $R_I=R-C(I)$. Then
\begin{equation}\label{eq:line-residual}
G[A\setminus N_G(I)]=L(R_I).
\end{equation}

\begin{lemma}\label{lem:residual-third}
For every $I\subseteq B$, the root $R_I$ supports the uniform fractional matching $x_e=1/3$.
\end{lemma}
\begin{proof}
For $U\subseteq V(R_I)$ the budget in Lemma~\ref{lem:third} is also
\begin{equation}\label{eq:original-budget}
P_{R_I}(U)+r_{R_I}(U)
=\sum_{T\in U}|T\cap B|+|\delta_R(U)|.
\end{equation}
Edges leading to deleted root vertices are counted in $P_{R_I}(U)$, which explains the equality. If an odd set violated~\eqref{eq:ports}, this sum would be one. Thus all triangles in $U$ consist of vertices represented by internal root edges, apart from exactly one attachment $z$: either a single $B$ occurrence or the $A$ vertex of the sole root edge leaving $U$.

Let $X$ be the original vertices represented by the internal root edges. It is nonempty, since $2e_R(U)=3|U|-1$ and a singleton cannot have budget one. Every neighbor of $X$ lies in $X\cup\{z\}$. There are vertices outside this set. If $z\in B$, it has only two neighbors in $X$, but has degree at least four. If $z\in A$, its triangle outside $U$ contains vertices outside $X\cup\{z\}$. Therefore $z$ is a cutvertex of $G$, a contradiction. Lemma~\ref{lem:third} applies.
\end{proof}

\begin{proof}[Proof of Theorem~\ref{thm:main-residual}]
Combine~\eqref{eq:line-residual}, Lemma~\ref{lem:residual-third}, and~\eqref{eq:uniform}; the independent set in the residual line graph can be adjoined to $I$.
\end{proof}

\begin{proposition}\label{prop:root-exact}
Let $\tau=|\mathcal T|$, and let $t(v)$ count triangles through $v\in B$. Then
\begin{align}
\alpha(G)&=\max_{I\in\IS(G[B])}\bigl(|I|+\nu(R_I)\bigr)\label{eq:root-exact}\\
&=\frac{\tau}{2}+\frac12\max_{I\in\IS(G[B])}
\left(\sum_{v\in I}(2-t(v))-\Def(R_I)\right).\label{eq:root-def}
\end{align}
Moreover $\Def(R_I)\le P_{R_I}(V(R_I))/3$.
\end{proposition}
\begin{proof}
Every independent set splits into its $B$ part $I$ and a matching of $R_I$, and conversely every such choice is compatible. No selected vertices share a triangle, so $|C(I)|=\sum_{v\in I}t(v)$. Substitute $2\nu(R_I)=\tau-|C(I)|-\Def(R_I)$.
\end{proof}

\section{A counterexample of minimum order}\label{sec:minimal}
Assume Conjecture~\ref{conj:main} is false. Choose a counterexample $G$ of minimum order and, subject to that condition, maximum size. Set
\[
n=|V(G)|,\qquad k=\alpha(G),\qquad d=5n-16k>0.
\]
Every $C_4$-free planar graph $H$ of order less than $n$ satisfies $\alpha(H)\ge\lceil5|V(H)|/16\rceil$. In particular, this holds for every proper induced subgraph of $G$. Minimum order is taken among all counterexamples, with no minimum-degree restriction. The maximum-size condition is used for the cofacial-edge assertion in Proposition~\ref{prop:chord}; the other reductions use vertex minimality.

\subsection{Deletion stability and separators}

\begin{lemma}\label{lem:stability}
For every $X$ with $1\le|X|\le3$, $\alpha(G-X)=k$. Moreover,
\[
\delta(G)\ge3,\qquad 1\le d\le5,\qquad k=\lfloor5n/16\rfloor.
\]
The residue $n\pmod{16}$ is respectively $13,10,7,4,1$ for $d=1,2,3,4,5$.
\end{lemma}
\begin{proof}
Minimality gives $\alpha(G-X)\ge k+(d-5|X|)/16>k-1$, so integrality proves stability. A vertex of degree at most two could be added to a maximum independent set of $G-N[v]$, a contradiction. Deleting one vertex gives $16k\ge5(n-1)$ and hence $d\le5$. The floor identity and residues follow from $5n=16k+d$.
\end{proof}
In particular, every prescribed triangle can be avoided by some maximum independent set. Thus no triangle is \emph{independence-associated}, in the sense of intersecting every maximum independent set.

\begin{theorem}\label{thm:cutsets}
$G$ is connected and has no clique cutset. It is $2$-connected, every triangle is facial in every plane embedding, and every facial boundary is a cycle.
\end{theorem}
\begin{proof}
Disconnectedness contradicts the sum of the bounds on smaller components. Suppose a clique $S$ of size $r\le3$ separates components with vertex sets $A_1,\ldots,A_p$. Put $a_i=\alpha(G[A_i])$ and $\sigma_i=16a_i-5|A_i|\ge0$. Deletion stability gives
\[
\sum_i a_i=k,\qquad \sum_i\sigma_i=5r-d<5r.
\]
Let $F_i=\{s\in S:\alpha(G[A_i\cup\{s\}])=a_i\}$. Since $S$ is a clique, $\alpha(G[A_i\cup F_i])=a_i$. Minimality gives $\sigma_i\ge5|F_i|$, so $\sum_i|F_i|<r$. Choose $s\in S\setminus\bigcup_iF_i$. For each $i$, a maximum independent set in $G[A_i\cup\{s\}]$ has size $a_i+1$ and contains $s$. Their union has size $k+1$, a contradiction. The facial assertions follow because a separating triangle is a clique cutset and $G$ is $2$-connected.
\end{proof}

\subsection{Neighborhood and boundary restrictions}
\begin{lemma}\label{lem:closed}
Every nonempty independent set $I$ satisfies
\begin{equation}\label{eq:closed}
5|N_G[I]|-16|I|\ge d.
\end{equation}
\end{lemma}
\begin{proof}
Extend an independent set of $G-N[I]$ by $I$ and use minimality:
$k\ge |I|+5(n-|N[I]|)/16$.
\end{proof}

\begin{proposition}\label{prop:incidence}
If $I$ consists of independent degree-three vertices, the bipartite incidence graph on $I$ and $N(I)$, retaining just the edges incident with $I$, is a forest. Every connected component contains at most $5-d$ vertices of $I$.
\end{proposition}
\begin{proof}
An incidence cycle containing $r$ selected vertices supplies an independent subset $J$ with at most $2r$ distinct neighbors, because each selected vertex has one remaining incident edge outside that cycle. Thus $|N[J]|\le3r$, contradicting~\eqref{eq:closed}. For a tree component with $r$ selected vertices there are $3r$ edges, hence $3r+1$ vertices. Apply~\eqref{eq:closed} to its selected vertices to obtain $d\le5-r$.
\end{proof}
This implies that a vertex has at most $5-d$ pairwise nonadjacent degree-three neighbors. If $d=5$, there are no degree-three vertices.

For an induced region of order $h$, external vertex boundary of size $b$, and an independent set of size $a$,~\eqref{eq:closed} gives
\begin{equation}\label{eq:region}
16a-5(h+b)\le-d.
\end{equation}
In particular, a sixteen-vertex induced region with boundary at most three has independence number exactly five: minimality gives the lower bound, and an independent set of size six would violate~\eqref{eq:region}.

\begin{proposition}\label{prop:chord}
Every nonadjacent cofacial pair of vertices is joined by a three-edge path in $G$. If $uvw$ is a triangle with $d_G(u)=d_G(v)=3$, their external neighbors $x,y$ are distinct and nonadjacent and are joined by a three-edge path in $G-\{u,v,w\}$.
\end{proposition}
\begin{proof}
Adding a cofacial edge preserves planarity and cannot increase independence number. Edge maximality forces a four-cycle containing the new edge, proving the first statement. For the second, Lemma~\ref{lem:neighbors} gives $x\ne y$ and $xy\notin E(G)$. Deleting the triangle and adding $xy$ preserves planarity. Every independent set of this modified graph extends by $u$ or $v$, since it cannot contain both $x,y$. If the modified graph were $C_4$-free, minimality would give $\alpha(G)\ge1+5(n-3)/16>5n/16$. Hence its new edge lies in a four-cycle, giving the path.
\end{proof}

\section{Planar applications in minimum degree four}\label{sec:applications}
We first derive a face-counting identity for a minimum counterexample with
$\delta(G)\ge4$, using the notation of Sections~\ref{sec:residual}
and~\ref{sec:minimal}. We then treat ordinary cycles and prove the
small-order exclusion in Theorem~\ref{thm:no-three-order}.

\subsection{Face counts and matching deficiencies}
Let $G$ be the minimum counterexample of Section~\ref{sec:minimal}, and assume $\delta(G)\ge4$. It is $2$-connected by Theorem~\ref{thm:cutsets}, so the root construction of Section~\ref{sec:residual} applies. Taking $I=\varnothing$ in Theorem~\ref{thm:main-residual} gives $\alpha(G)\ge\lceil|A|/3\rceil$. Since $\alpha(G)<5n/16$, we obtain
\begin{equation}\label{eq:many-claws}
|B|>n/16.
\end{equation}
There is also a useful exact face identity. Write $a=|A|$, $b=|B|$,
\[
s=\sum_{f:\ell(f)\ne3}(\ell(f)-5),\qquad D=\sum_v(d_G(v)-4).
\]
By Theorem~\ref{thm:cutsets}, all triangles are facial and every face boundary is a cycle. Thus $s,D\ge0$. If $m=|E(G)|$ and $\tau$ is the number of triangles, degree and triangle-incidence counts give
\[
2m=4n+D,\qquad 3\tau=2a+\sum_{v\in B}t(v).
\]
Euler's formula and the sum of face lengths give
$3m-2\tau=5n-10-s$. Substituting the preceding identities yields
\begin{equation}\label{eq:charge-zero}
a-3b+2\sum_{v\in B}t(v)=30+3s+\frac92D.
\end{equation}
For any proper four-coloring $B=I_1\cup\cdots\cup I_4$, put $\Delta_i=\Def(R_{I_i})$. Averaging~\eqref{eq:root-def} and using~\eqref{eq:charge-zero} yields
\begin{equation}\label{eq:four-def}
\alpha(G)\ge\frac{5n}{16}+\frac58+\frac{s}{16}+\frac{3D}{32}
-\frac18\sum_{i=1}^4\Delta_i.
\end{equation}
Thus it would suffice to find a coloring with $\sum_i\Delta_i\le5+s/2+3D/4$. The uniform-third lemma does not by itself prove this stronger estimate.

\subsection{Ordinary cycles}
Call an edge \emph{ordinary} if it belongs to no triangle. The next theorem is an unconditional application of Theorem~\ref{thm:main-residual}; in its statement, $G$ need not be a minimum counterexample.
\begin{theorem}\label{thm:c5}
Suppose $G$ is $2$-connected, $C_4$-free, planar, and $4$-regular, and its ordinary edges form exactly one five-cycle. Then
\[
\alpha(G)\ge\frac{|V(G)|-1}{3}\ge\frac{5|V(G)|}{16}.
\]
\end{theorem}
\begin{proof}
The five cycle vertices form $B$, each in one triangle, and every other vertex is in two triangles. The five triangles are distinct: adjacent cycle vertices cannot share a triangle because their edge is ordinary; two vertices at distance two sharing a triangle would make an ordinary edge part of another triangle. Thus $G[B]=C_5$ and $\tau=(2n-5)/3$, giving $n\equiv1\pmod3$.

There are $n+2$ faces. The number of triangular faces is at most $\tau$, and all other faces have length at least five, so $4n\ge5(n+2)-2\tau$, giving $n\ge40$. Choose a nonadjacent pair $I\subseteq B$. Its vertices have at most four $A$ neighbors altogether. Theorem~\ref{thm:main-residual} gives
\[
\alpha(G)\ge\left\lceil2+\frac{n-9}{3}\right\rceil=\frac{n-1}{3}.
\]
Finally $(n-1)/3-5n/16=(n-16)/48\ge0$.
\end{proof}
For the analogous ordinary-six-cycle family, its six vertices lie in distinct triangles. Opposite vertices cannot share a triangle because a three-edge cycle arc would give a four-cycle. Selecting an alternating triple gives $\alpha(G)\ge(n-3)/3$. Euler counting gives $n\ge42$ and triangle incidence gives $n\equiv0\pmod3$. Hence this estimate proves~\eqref{eq:target} for that family when $n\ge48$. Theorem~\ref{thm:no-three-order} treats the remaining orders forty-two and forty-five under the additional assumption that $G$ is a minimum counterexample.

\subsection{Small roots with cofacial ports}
We first prove two matching lemmas independent of the minimum-counterexample
assumption. A degree-two vertex in these roots will be called a \emph{port}.
The lemmas give small-order information not supplied by the uniform-third
bound alone.

\begin{lemma}\label{lem:cofacial-ports}
Let $H$ be a connected simple planar graph of girth at least five, with all
degrees two or three. If at most four vertices have degree two, and these
vertices are incident with a common face, then $|V(H)|\ge15$.
\end{lemma}
\begin{proof}
Write $u=|V(H)|$ and let $p\le4$ count its degree-two vertices. Suppose $u\le14$.
A bridge has two cyclic sides: a tree side would have internal degree deficit
$p_i+1$, giving $u_i=p_i-1$, impossible since $p_i\le u_i$.
The girth bound on both sides gives $u_i\ge17-3p_i$, hence $u\ge34-3p\ge22$.
Thus $H$ is bridgeless and, being subcubic, $2$-connected.

Choose the common face as outer, with boundary length $L$. Euler counting gives
\[
u\ge20-3p,\qquad L\le\frac{u+3p-10}{2}.
\]
Let $r=u-L$ count the vertices strictly inside the outer boundary, let $c$
count boundary chords, and let $h$ count edges joining boundary vertices to
interior vertices. Since all $p$ ports are on the boundary, its degree sum is
$3L-p=2L+2c+h$, so $h=L-p-2c$. Every interior vertex has degree three.
Consequently, the number $e_{\rm int}$ of edges in the subgraph induced by
these $r$ interior vertices satisfies
\[
2e_{\rm int}=3r-h=3u+p-4L+2c,
\qquad e_{\rm int}\ge\frac{3u+p-4L}{2}.
\]
For $r\ge4$, the planar girth bound $e_{\rm int}\le5(r-2)/3$ yields
$u+2L\ge20+3p$. But the displayed upper bound on $L$ and $u\le14$ give
$u+2L\le18+3p$, a contradiction. If $r\le3$, then $L\ge u-3$ implies
$u\le3p-4$. Together with $u\ge20-3p$, $p\le4$, and $u\equiv p\pmod2$,
the only possibility is $u=8$, $p=4$, $L=5$, $r=3$. It gives
$e_{\rm int}\ge4$, impossible on three vertices in a girth-five graph.
\end{proof}

\begin{lemma}\label{lem:small-perfect}
Let $H$ be a connected simple planar graph of girth at least five, with all
degrees two or three and at most four degree-two vertices. If its order is
even and at most thirty-two, then it has a perfect matching.
\end{lemma}
\begin{proof}
By Tutte's theorem~\cite{Tutte1947}, failure gives a nonempty set $S$ of order
$s$ such that $H-S$ has $q\ge s+2$ odd components. For an odd component $W$,
let $p_W$ count its original degree-two vertices and $r_W$ its edges to $S$.
Connectivity gives $r_W\ge1$, and
\[
2e(W)=3|W|-p_W-r_W,\qquad p_W+r_W\text{ is odd},\qquad
\sum_W r_W\le3s.
\]
If $W$ is cyclic, its internal degree deficit is $p_W+r_W$, so the planar
girth bound gives
\[
|W|\ge20-3(p_W+r_W).
\]
If $W$ is a tree, the degree identity instead gives
$|W|=p_W+r_W-2$.
Let $q_1,q_2$ count the odd components with respectively one and two edges to
$S$. The cut-capacity inequalities
$3q-2q_1-q_2\le\sum_Wr_W\le3s$, together with $q\ge s+2$,
give $2q_1+q_2\ge6$. Each cut-two component
has odd $p_W$, so $q_2\le4$. A cut-one component is cyclic: the tree degree
identity would give $|W|=p_W-1$. Its order is at least $17-3p_W$, and $p_W$
is even.

If $q_1\ge3$, the cut-one components alone total at least $51-12=39$ vertices.
If $q_1=2$ and $q_2\ge3$, all cut-one port counts are zero, so those two
components total at least thirty-four. If $q_1=2,q_2=2$, the cut-one components
total at least twenty-eight. Cut capacity is then tight: $q=s+2$, every other
odd component has three exits, and $s\ge2$. Including the cut-two components,
the other odd components, and $S$ gives at least $28+2s\ge32$ vertices.
Equality forces $s=2$ and two singleton cut-two components. Both singletons
have the same two neighbors in $S$, giving a four-cycle.

The remaining possibility is $q_1=1,q_2=4$. All four degree-two vertices lie
in those cut-two components; the cut-one component has at least seventeen
vertices. Again capacity is tight: $q=s+2$, all other odd components have
three exits, all vertices of $S$ have degree three, $S$ is independent, and
there are no even components. Here $s\ge3$. A nonsingleton component with
three exits and no ports is cyclic and has order at least eleven; it would
force total order at least $2s+28\ge34$. A cut-two nonsingleton with one port
has the same lower bound. Hence all components except the cut-one component
are singletons. The cut-one component is joined to the rest of $H$ by a
single edge, which is a bridge. Since $H$ is connected, deleting that
component leaves the other connected side of this bridge. The remaining
graph is bipartite, with parts $S$ and the singleton components of $H-S$,
and is planar of girth at least six. It has $2s+1$ vertices and $3s-1$ edges. The planar
girth-six bound gives $3s-1\le3(2s-1)/2$, a contradiction.
\end{proof}

\subsection{The order of a minimum counterexample}
We now return to the minimum counterexample of Section~\ref{sec:minimal}.
The following statement uses minimum order among all counterexamples,
not only among graphs of minimum degree at least four.

\paragraph{Embedded roots.}
Fix a plane embedding of $G$. All its triangles are facial by
Theorem~\ref{thm:cutsets}. In each triangular face place a center and draw
three noncrossing arms to its corners. Centers in different faces and
interiors of their arms are disjoint. To form a root, retain the arms
corresponding to represented vertices; the two arms ending at such a
vertex are joined there to make a root edge. Suppressing this degree-two
junction gives the root embedding. The original edges of $G$ are not
part of this root drawing. A triangle with one removed corner retains its
center and two arms, and therefore becomes a degree-two node representing
the remaining edge clique. Its omitted arm is still available as an arc
from that port to the removed corner, disjoint from the root except at
the port. A triangle whose corners are all removed contributes no root
node. This construction realizes the clique root in the original
embedding, rather than choosing an unrelated embedding of its abstract
graph. We apply Lemma~\ref{lem:cofacial-witness} to the omitted arms below.

\begin{theorem}\label{thm:no-three-order}
If $G$ is a counterexample of minimum order to Conjecture~\ref{conj:main}
and has no degree-three vertex, then $|V(G)|\ge49$.
\end{theorem}
\begin{proof}
By Lemma~\ref{lem:stability}, the degree assumption implies $\delta(G)\ge4$.
Use the root notation of Section~\ref{sec:residual} and the face-counting
parameters defined in Section~\ref{sec:applications}.
From~\eqref{eq:charge-zero},
\[
\begin{aligned}
n&=30+3s+\sum_{v\in B}\left(4-2t(v)+\frac92(d_G(v)-4)\right)\\
 &\ge30+3s+2b+\frac52D.
\end{aligned}
\]
For degree four, a claw center has $t(v)\le1$; for degree at least five,
$t(v)\le\lfloor d_G(v)/2\rfloor\le d_G(v)-3$, proving the inequality.
The integer $D$ is even. If $D>0$, then $D\ge2$. Since $s\ge0$ and
$b>n/16$ by~\eqref{eq:many-claws}, the coefficient $5/2$ of $D$ gives
\[
n\ge30+2b+\frac52D\ge35+2b>35+\frac n8.
\]
Thus $7n/8>35$, or $n>40$.
If $D=0$, the ordinary-edge graph on $B$ has degrees two or four and girth
at least five, so $b\ge5$ and $n\ge40$. The residue condition excludes
orders forty and forty-one. Below forty-nine, only forty-two and forty-five
remain.

The degree budget at those two orders can be classified without graph
enumeration. We have $b\ge3$. If $D\ge4$, the last inequality gives $n\ge46$.
For $D=2$, the degrees above four are either two fives or one six. Write $r$
for the total loss of triangle incidences from their maximum values, also
counting degree-four claw centers with no triangle. The charge identity is
\[
n=35+2b+2r+3s.
\]
It has no solution at $n=42$ with $b\ge3$. At $n=45$, it requires $s=0$ and
$b+r=5$. The ordinary-edge graph has four edges and girth at least five,
hence is a forest. This forces $b=5,r=0$ and the two degree-five vertices
are the ends of an ordinary five-vertex path. The one-degree-six alternative
would give an isolated ordinary-edge vertex in a five-vertex tree and is
impossible.

For $D=0$, put $r=|\{v\in B:t(v)=0\}|$, so
$n=30+2b+2r+3s$. The ordinary-edge graph has minimum degree two, giving
$5\le b\le7$ in the range considered. On five vertices it must be a five-cycle:
any additional edge is a forbidden chord. Direct substitution leaves only
$b=6,r=0$, with $s=0$ at order forty-two and $s=1$ at order forty-five.
Thus the ordinary edges form a single six-cycle.

\emph{The six-cycle cases.}
The six cycle vertices lie in distinct triangles. Its triangle root has
$(2n-6)/3$ vertices, six degree-two ports, and all other degrees three.
Use the root embedding just described, retaining the arms at the $A$
vertices. The drawing of the ordinary cycle is disjoint from this root:
its vertices are not represented by root edges, and its edges avoid the
interiors of all triangular faces. Each of the
six ports is joined to its cycle vertex by the omitted arm of that
triangle. The interior of this arm misses the entire root. Thus the
ordinary cycle is a connected witness $J$ in
Lemma~\ref{lem:cofacial-witness}; all ports are incident with a common face
of the root, and the ports of each component are cofacial in its inherited
embedding.

A component with no attachment would give a component of $G$ disjoint
from the ordinary cycle. A component with only one attachment would make
that cycle vertex a cutvertex of $G$. Hence every root component has at
least two ports. If there are multiple components, each has two to four
ports, since there are six altogether. Applying
Lemma~\ref{lem:cofacial-ports} to each component forces at least thirty
root vertices, more than twenty-six or twenty-eight. Hence the root is
connected. To see that the ordinary cycle is facial, use the same embedding,
whose root edges do not cross the ordinary cycle.
Each of its six vertices belongs to exactly one triangle, lying on one side
of the cycle. If the cycle were separating, vertices on both sides would
belong to triangle-root components on their respective sides, with no root
edge joining them. This would contradict root connectivity. The cycle has
no chord, since any chord of a six-cycle creates a triangle containing an
ordinary edge or a four-cycle. Its empty side is therefore a face.
At order forty-two all nontriangular faces are pentagons, so a facial
six-cycle is impossible.

At order forty-five, replace each triangle and the facial six-cycle by a face
center and represent every original vertex by the edge between its two selected
faces. This gives a simple plane graph $X$ on twenty-nine vertices and
forty-five edges, with one degree-six vertex $z$, all other degrees three,
and all faces pentagonal. The graph $G$ is its medial graph. No bridge of
$X$ is possible, because its representing vertex would be a cutvertex of
$G$. A degree-three cutvertex of $X$ would require a bridge, and $z$ is not
a cutvertex because $X-z$ is the connected triangle root. Thus $X$ is
$2$-connected.

There is no two-edge bond in $X$. For a connected side $J$ avoiding $z$,
the degree identity is $e(J)=(3|J|-2)/2$. This side is cyclic and lies in the
girth-five triangle root, so $|J|\ge14$. The other side $K$, containing $z$,
has $e(K)=(3|K|+1)/2$. All its faces except the one containing the deleted
connected region $J$ are pentagons. Choose that face as outer, of boundary
length $L\ge3$. Euler gives $3e(K)=5|K|-5-L$, implying $|K|\ge19$.
Together these force $|V(X)|\ge33$, impossible. In plane duality, bridges
correspond to loops and two-edge bonds correspond to pairs of parallel
dual edges. Thus the dual $X^*$ has neither loops nor parallel edges.
It is connected because $X$ is connected, and it is $2$-connected because
plane duality preserves $2$-connectedness (equivalently, a nontrivial
block decomposition on either side gives one on the other).
Since the faces of $X$ are pentagons, $X^*$ is $5$-regular; its faces
correspond to the vertices of $X$, giving twenty-eight triangular faces
and one hexagonal face. Keith, Froncek, and Kreher's nonexistence
result~\cite{Keith2018}, as stated in their corrigendum~\cite[Theorem~1]{KeithCorr2018},
excludes a finite planar regular $2$-connected graph with all but one face
of the same length and a single face of a different length. The dual just
constructed satisfies every hypothesis, giving a contradiction.

\emph{The ordinary-path case at order forty-five.}
Write the path as $x b_1 b_2 b_3 y$, with degree-five ends in two triangles
each, degree-four internal vertices in one triangle each, and every other
vertex in two triangles. Let $\{b_2,a,c\}$ be the middle triangle and delete
$N[b_2]=\{b_1,b_2,b_3,a,c\}$. A retained vertex loses at most one neighbor,
since otherwise it and $b_2$ have two common neighbors. The remaining graph
$K$ is claw-free, with degrees three or four: the ordinary edge at each of
$x,y$ has been deleted, and every other retained vertex was already a
noncenter in $G$.

Triangle incidence gives $3\tau=2\cdot45-3=87$, so $G$ has twenty-nine
triangles. The affected triangles are the unique triangles through
$b_1,b_2,b_3$ and the second triangles through $a$ and $c$. They are five
distinct triangles. Indeed, the only edge in $N(b_2)$ is $ac$, by
Lemma~\ref{lem:neighbors}; thus no other triangle contains two deleted
neighbors of $b_2$, and the second triangles through $a,c$ cannot coincide.
The triangle $\{b_2,a,c\}$ disappears completely, whereas each of the other
four loses exactly one vertex and becomes a maximal edge clique. The other
twenty-four triangles are unchanged.

Every retained vertex belongs to two surviving maximal cliques, counting
these four edge cliques. The clique root $H$ consequently has twenty-eight
vertices: twenty-four triangle nodes of degree three and four edge-clique
nodes of degree two. There are no private leaves. Use the embedded-root
construction above: keep the centers of the twenty-four intact triangles
and of the four triangles with one deleted corner, with arms only to
retained vertices. This is the clique root $H$ of $K$. The original drawing
of $J=G[N[b_2]]$ is connected, because it contains every edge from $b_2$
to its four neighbors, and is disjoint from the drawing of $H$. Each of
the four ports has an omitted arm to its deleted corner, which is a
vertex of $J$; the arm meets $H$ only at the port. Applying
Lemma~\ref{lem:cofacial-witness} with this connected witness $J$ proves
that all four ports are incident with a common face of $H$, including
when $H$ is disconnected. In particular, the ports of every component
are cofacial in the inherited embedding.

The root is simple, planar, and of girth at least five by
Lemma~\ref{lem:line-root}. If it were disconnected, each component would
satisfy Lemma~\ref{lem:cofacial-ports}, with at most four ports; for a
component with no ports the face condition is vacuous. There would be
at least thirty vertices, a contradiction.
Thus $H$ is connected, and Lemma~\ref{lem:small-perfect} gives a perfect
matching of size fourteen. Adding $b_2$ to the corresponding independent set
of $K$ gives $\alpha(G)\ge15$, contradicting $k=\lfloor225/16\rfloor=14$.

Both order-forty-five types are excluded. The residue restriction then excludes
forty-six, forty-seven, and forty-eight, completing the proof.
\end{proof}

\section{Triangle tree components and their boundary states}\label{sec:branches}
Throughout this section, $G$ is again the minimum counterexample of Section~\ref{sec:minimal}; degree-three vertices are now allowed. We determine how a class of degree-three components interacts with independent sets in its complement.

\subsection{Exact boundary states}
A \emph{triangle tree component} is a component $T$ of $G[V_3]$, where $V_3=\{v:d_G(v)=3\}$, consisting of $q$ disjoint triangles joined by $q-1$ bridges whose contraction gives a tree. Every bridge uses distinct triangle vertices, since all vertices have degree three. There are $q+2$ external ports. Put $F=G-V(T)$ and $U=N_G(V(T))\setminus V(T)$; external neighbor identifications are allowed.

\begin{lemma}\label{lem:profile}
For every $J\in\IS(F)$,
\begin{equation}\label{eq:profile}
\alpha(T-N_G(J))=q-1_{\{U\subseteq J\}}.
\end{equation}
Moreover $q\le5-d$.
\end{lemma}
\begin{proof}
If some external port is available, root the contracted tree at that port's triangle and select the port. In each other triangle select the endpoint of its bridge to the parent. This gives one selected vertex per triangle and no selected adjacent pair. If every external port is blocked, selecting $q$ vertices would select $q$ bridge endpoints, whereas only $q-1$ bridges exist and at most one endpoint of each can be selected. The same parent-facing construction in the nonroot triangles gives $q-1$.

For the bound on $q$, select an external root port and all parent-facing ports. Its closed neighborhood consists of all $3q$ vertices of $T$ and the single external neighbor of the root port. Therefore~\eqref{eq:closed} gives $d\le5(3q+1)-16q=5-q$.
\end{proof}

\begin{theorem}\label{thm:forced}
The terminal set $U$ is independent, and
\begin{equation}\label{eq:forced}
\alpha(F)=k-q+1.
\end{equation}
Every maximum independent set of $F$ contains all of $U$. For $1\le|X|\le3$ with $X\cap U\ne\varnothing$,
\begin{equation}\label{eq:robust}
\alpha(F-X)=\alpha(F)-1.
\end{equation}
Every distinct pair in $U$ is joined by a three-edge path in $F$.
\end{theorem}
\begin{proof}
The profile gives $\alpha(F)\le k-q+1$. Minimality gives
\[
\alpha(F)\ge k-q+\left\lceil\frac{d+q}{16}\right\rceil=k-q+1,
\]
since $1\le d+q\le5$. A maximum set omitting a terminal would extend by $q$, contradicting $\alpha(G)=k$. Thus all terminals are forced and are independent. For~\eqref{eq:robust}, forced membership gives the upper bound, while minimality gives
\[
\alpha(F-X)\ge k-q+\left\lceil\frac{d+q-5|X|}{16}\right\rceil=k-q.
\]

After deleting the connected region $T$, its terminals are incident with a common face. A new edge between two distinct terminals can be drawn through that region. Since they are forced in every maximum independent set of $F$, adding this edge reduces its independence number by exactly one. If the resulting graph were $C_4$-free, minimality at order $n-3q$ would contradict~\eqref{eq:forced}. Hence the new edge lies in a four-cycle, supplying a three-edge path in $F$. This does not assert that the graph distance is three.
\end{proof}

\subsection{Restrictions on terminals and exchanges}
The forced membership in Theorem~\ref{thm:forced} imposes both local degree restrictions and compatibility conditions between different components.

\begin{corollary}[Degrees around terminals]\label{cor:terminal-degrees}
For every $v\in U$, deleting any set of at most two vertices from $F-v$ preserves its independence number. Consequently,
\[
\delta(F)\ge2,\qquad d_F(w)\ge3\quad\text{for every }w\in N_F(U).
\]
If $h(u)$ external ports of $T$ have the same neighbor $u\in U$, then $d_G(u)\ge h(u)+2$.
\end{corollary}
\begin{proof}
The first assertion follows from~\eqref{eq:robust} by deleting $v$ together with at most two other vertices. A graph whose independence number survives every deletion of at most two vertices has minimum degree at least two: otherwise delete the closed neighborhood of a vertex of degree at most one and then add that vertex back to an independent set of unchanged size. Thus $\delta(F-v)\ge2$ for each $v\in U$. The set $U$ contains at least two distinct vertices, since the two external ports of a leaf triangle have distinct neighbors by $C_4$-freeness. For any $w\in V(F)$, choose $v\in U\setminus\{w\}$ to obtain $d_F(w)\ge2$. If $w$ is adjacent to a terminal $v$, the inequality $d_{F-v}(w)\ge2$ gives $d_F(w)\ge3$. Finally $d_G(u)=d_F(u)+h(u)$.
\end{proof}

\begin{proposition}[Incomparable terminal sets]\label{prop:antichain}
For distinct triangle tree components $T_j,T_l$, their terminal sets satisfy
\[
U_j\nsubseteq U_l,\qquad U_l\nsubseteq U_j.
\]
More precisely, if $U_j\cap U_l\ne\varnothing$, no maximum independent set of $G-T_j$ contains all of $U_l$.
\end{proposition}
\begin{proof}
Every maximum independent set $M$ of $G-T_j$ contains $U_j$. Suppose it also contains $U_l$, and choose $u\in U_j\cap U_l$. The profile of $T_l$ gives $|M\cap V(T_l)|=q_l-1$: otherwise that part could be enlarged, contradicting maximal size. Remove $u$ and replace this part by $q_l$ compatible vertices of $T_l$, using the newly available port at $u$. The resulting independent set has size $|M|$ and omits the forced terminal $u$, a contradiction. If $U_l\subseteq U_j$, every maximum set $M$ would contain $U_l$, so this containment is impossible. Interchange $j,l$ to exclude the other containment.
\end{proof}

\begin{corollary}\label{cor:branch-claw}
For every triangle tree component $T$, its complement $F$ contains an induced claw. More generally, if deleting $c$ vertices of $F$ makes it claw-free, then
\begin{equation}\label{eq:branch-claw}
5n+16q+31d\le160(c+1).
\end{equation}
\end{corollary}
\begin{proof}
Choose $v\in U$, so $\alpha(F-v)=k-q$. If $F-S$ is claw-free and $|S|=c$, then $F-(S\cup\{v\})$ is claw-free and has at least $n-3q-c-1$ vertices. Theorem~\ref{thm:main-claw} gives
\[
k\ge q+\left\lceil\frac{10(n-3q-c-1)}{31}\right\rceil
=\left\lceil\frac{10n+q-10c-10}{31}\right\rceil.
\]
Substitute $16k=5n-d$ to obtain~\eqref{eq:branch-claw}. If $c=0$ and $n\ge29$, then $5n+16q\ge160$ contradicts the strict target failure. For $n=26,27,28$ the displayed lower bound gives $k\ge9$, whereas $k=\lfloor5n/16\rfloor=8$. The Ramsey reduction excludes smaller orders.
\end{proof}

\begin{lemma}[Exchange connectivity]\label{lem:exchange}
Fix a maximum independent set $M$ of $F$. For any nonempty $X\subseteq U$ with $|X|\le3$, there exist disjoint sets $A_X\subseteq M$ and $B_X\subseteq V(F)\setminus M$ such that
\[
X\subseteq A_X,\quad |A_X|=|B_X|+1,
\quad (M\setminus A_X)\cup B_X\in\IS(F),
\]
and $F[A_X\cup B_X]$ is connected and bipartite. For each $u\in U\cap A_X$, that bipartite graph has a matching saturating $B_X$ and avoiding $u$.
\end{lemma}
\begin{proof}
By~\eqref{eq:robust}, choose an independent set $J$ of size $|M|-1$ avoiding $X$. Each component of $F[M\triangle J]$ is bipartite, with parts in $M\setminus J$ and $J\setminus M$. Its first part has size at least its second, since otherwise swapping that component increases $M$. If a component contains a forced vertex lost from $M$, equality is impossible: the swap would produce another maximum set omitting that vertex. The total imbalance is one, so all lost forced vertices lie in a single connected component with imbalance one. Swap only this component to define $A_X,B_X$.

For $B'\subseteq B_X$, maximality of $M$ gives $|N_F(B')\cap M|\ge|B'|$. If this neighborhood contains a forced vertex, equality would again give a maximum set omitting it, so the inequality is strict. All these neighbors lie in $A_X$. Deleting any $u\in U\cap A_X$ therefore preserves Hall's inequalities for matching $B_X$ into $A_X\setminus\{u\}$. Hall's theorem proves the claim.
\end{proof}
This exchange restriction supplements the three-edge-path condition. The latter alone does not encode the robustness in~\eqref{eq:robust}.

\section{A joint model for branches and claw centers}\label{sec:joint}
Continue with the minimum counterexample of Section~\ref{sec:minimal}, and assume in addition that every component of $G[V_3]$ is a triangle tree $T_j$ with $q_j$ triangles. Put $Q=\sum_j q_j$, delete all these components, and let $F$ be the remaining core. Every core vertex has original degree at least four. Partition it into the original claw centers $B$ and the other vertices $A$.

\subsection{The exact matching formulation}
\begin{lemma}\label{lem:terminals-B}
Every terminal set $U_j=N_G(T_j)\setminus V(T_j)$ is contained in $B$.
\end{lemma}
\begin{proof}
An external edge of a triangle tree belongs to no triangle. Otherwise its degree-three endpoint would share a second triangle with one of its two triangle edges, contradicting $C_4$-freeness. Its other endpoint has degree at least four, because $T_j$ is a whole component of $G[V_3]$. A non-claw-center of degree at least four has degree four and all incident edges triangular. Hence that endpoint lies in $B$.
\end{proof}
Let $R$ have a node for every core triangle and an edge for every $A$ vertex. The core triangles are unaffected by deleting the triangle trees. For an independent $I\subseteq B$, define $C(I)$ and $R_I$ as before. Set
\[
\mathcal F=\{I\in\IS(F[B]):U_j\nsubseteq I\text{ for every }j\}.
\]

\begin{theorem}\label{thm:joint}
Under these assumptions,
\begin{align}
\alpha(G)&=Q+\max_{I\in\mathcal F}\bigl(|I|+\nu(R_I)\bigr)\label{eq:joint}\\
&=Q+\frac{\tau}{2}+\frac12\max_{I\in\mathcal F}
\left(\sum_{v\in I}(2-t(v))-\Def(R_I)\right),\label{eq:joint-def}
\end{align}
where $\tau$ is the number of core triangles and $t(v)$ counts core triangles through $v$. Every $R_I$ supports the uniform fractional matching $1/3$. Consequently, with $a=|A|$ and $b=|B|$,
\begin{equation}\label{eq:joint-third}
\alpha(G)\ge Q+\max_{I\in\mathcal F}
\left(|I|+\frac{|A\setminus N_F(I)|}{3}\right),
\qquad b>n/16.
\end{equation}
\end{theorem}
\begin{proof}
For $J\in\IS(F)$, the profiles give an extension of size
\[
|J|+Q-\#\{j:U_j\subseteq J\}.
\]
The maximum of the penalized core objective equals its maximum subject to avoiding every full terminal set. Indeed, remove one terminal vertex for each violated constraint. At most the number of violations is removed, and removing vertices creates no new violation. Conversely every feasible set incurs zero penalty. Since $U_j\subseteq B$, feasibility depends only on $I=J\cap B$. The line-graph translation then proves~\eqref{eq:joint}, and the deficiency substitution proves~\eqref{eq:joint-def}.

The proof of Lemma~\ref{lem:residual-third} applies in the original graph $G$, even if the core is not $2$-connected: vertices represented by internal root edges are still degree-four vertices with two triangles and have no neighbors in deleted branches. A budget-one region attaches only through one original $A$ or $B$ vertex. That vertex has original degree at least four, so it is a cutvertex of $G$, impossible. Thus the uniform-third bound holds. The empty set belongs to $\mathcal F$, giving $\alpha(G)\ge Q+a/3=(n-b)/3$, and the strict counterexample inequality implies $b>n/16$.
\end{proof}

\subsection{Critical constraints and terminal restoration}
The boundary conditions in~\eqref{eq:joint} are individually critical: relaxing any one of them increases the optimum by exactly one. This gives a further restriction on the residual matchings.

\begin{proposition}[Critical boundary constraints]\label{prop:critical}
Write $f(I)=|I|+\nu(R_I)$ and $K=k-Q$. Then $\max_{I\in\mathcal F}f(I)=K$. For each branch $j$,
\begin{equation}\label{eq:critical}
\max\{f(I):I\in\IS(F[B]),\ U_j\subseteq I,\ U_l\nsubseteq I\ (l\ne j)\}=K+1.
\end{equation}
Thus removing any one constraint $U_j\nsubseteq I$, while retaining all others, increases the optimum by exactly one.
\end{proposition}
\begin{proof}
Choose a maximum independent set $M$ of $G-T_j$. Its size is $k-q_j+1$, and it contains all of $U_j$. For another branch $l$ whose terminal set is fully selected, Proposition~\ref{prop:antichain} gives $U_l\cap U_j=\varnothing$. Delete one selected terminal of $U_l$ and replace the $q_l-1$ selected branch vertices by $q_l$ compatible vertices. This preserves size and all of $U_j$, and destroys that full terminal selection. Deleting a core vertex creates no new full selection. Repeat until every remaining branch has $q_l$ selected vertices and no terminal constraint other than $j$ is violated. Deleting those branch vertices leaves an independent core set of size $k-Q+1$, containing $U_j$ and satisfying every other constraint. Its $B$ part supplies the lower bound in~\eqref{eq:critical}. For the upper bound, any such core choice extends to $G$ by $Q-1$ vertices, so its size is at most $k-Q+1$. The first equality follows from Theorem~\ref{thm:joint}.
\end{proof}

\begin{corollary}[Matching neutrality of terminal restoration]\label{cor:neutral}
Let $I_j$ attain the maximum in~\eqref{eq:critical}. For every $v\in U_j$,
\begin{align}
f(I_j\setminus\{v\})&=K,\qquad
\nu(R_{I_j\setminus\{v\}})=\nu(R_{I_j}),\label{eq:neutral}\\
\Def(R_{I_j\setminus\{v\}})&=\Def(R_{I_j})+t(v).\label{eq:neutral-def}
\end{align}
Writing $c_v$ for the number of root edges restored by removing $v$ from
$I_j$, we also have
\begin{equation}\label{eq:neutral-slack}
3\nu(R_{I_j})-|E(R_{I_j})|\ge c_v.
\end{equation}
Every maximum matching of $R_{I_j}$ remains maximum after that restoration
and leaves all the restored root vertices unmatched.
More generally, for every nonempty $X\subseteq U_j$,
\begin{equation}\label{eq:restore-set}
0\le\nu(R_{I_j\setminus X})-\nu(R_{I_j})\le |X|-1.
\end{equation}
\end{corollary}
\begin{proof}
The set $I_j\setminus\{v\}$ satisfies every boundary constraint. Its root
contains $R_{I_j}$, so its matching number cannot decrease. Thus its objective
is at least $f(I_j)-1=K$, and feasibility makes it at most $K$. Equality
forces the matching numbers to coincide. Independence of $I_j$ makes the
$t(v)$ restored triangle nodes disjoint from the nodes deleted by its other
vertices, giving~\eqref{eq:neutral-def}. Apply the uniform-third bound to the
restored root to obtain~\eqref{eq:neutral-slack}. An old maximum matching
uses none of the restored vertices and has the new maximum size.
For a nonempty set $X$, feasibility of $I_j\setminus X$ gives
$|I_j|-|X|+\nu(R_{I_j\setminus X})\le K=|I_j|+\nu(R_{I_j})-1$,
proving~\eqref{eq:restore-set}; root inclusion gives its lower bound.
\end{proof}

\subsection{Feasible colorings and the deficiency budget}
Each terminal set contains at least two vertices: the two external ports of a leaf triangle have distinct neighbors by Lemma~\ref{lem:neighbors}. Choose a pair in each $U_j$, and join that pair by a path whose internal vertices lie in the connected component $T_j$, using its two terminal edges. Distinct triangle-tree components are vertex-disjoint, so the interiors of these paths are disjoint in the fixed plane embedding. Paths may share a terminal, but they have no other common vertices or crossings. After deleting all branches, draw one virtual edge along each path; at a common terminal the edges retain their original incident sectors. The virtual edges can therefore be added simultaneously without crossings, giving a planar supergraph of the core. Repeated virtual edges may be suppressed. A proper four-coloring of its subgraph on $B$ partitions $B$ into four members of $\mathcal F$.

\begin{proposition}\label{prop:joint-charge}
Let $D_B=\sum_{v\in B}(d_G(v)-4)$ and let $s$ be the sum of the excesses over five of the nontriangular faces of $G$. Then
\begin{equation}\label{eq:joint-charge}
a-3b+2\sum_{v\in B}t(v)=30+3s+\frac92D_B-\frac{21}{2}Q.
\end{equation}
For any four-coloring respecting the virtual edges, set $\Delta_i=\Def(R_{I_i})$. Then
\begin{equation}\label{eq:joint-four}
\alpha(G)\ge\frac{5n}{16}+\frac58+\frac{s}{16}
+\frac{3D_B}{32}-\frac{5Q}{32}-\frac18\sum_i\Delta_i.
\end{equation}
Thus a counterexample must satisfy, for every such coloring,
\begin{equation}\label{eq:missing}
\sum_i\Delta_i>5+\frac{s}{2}+\frac{3D_B}{4}-\frac{5Q}{4}.
\end{equation}
\end{proposition}
\begin{proof}
There are $n=a+b+3Q$ vertices, $m=(4a+4b+9Q+D_B)/2$ edges, and $\tau+Q$ triangles, where $3\tau=2a+\sum_Bt(v)$. All triangles are facial. The face identity is
\[
3m-2(\tau+Q)=5n-10-s.
\]
Substitution gives~\eqref{eq:joint-charge}. Average the four feasible values in~\eqref{eq:joint-def}:
\[
\alpha(G)\ge Q+\frac{\tau}{2}
+\frac18\sum_{v\in B}(2-t(v))-\frac18\sum_i\Delta_i.
\]
Use the charge identity to obtain~\eqref{eq:joint-four}. Strict failure of~\eqref{eq:target} implies~\eqref{eq:missing}.
\end{proof}

\section{Discussion and further directions}\label{sec:next}
The matching approach separates two issues. Theorem~\ref{thm:main-claw} determines the sharp weighted bound when there are no claw centers, while Theorem~\ref{thm:main-residual} controls the residual graph after a compatible selection of centers. The remaining difficulty is to choose those centers so that the gains from the selection compensate for the matching deficiencies and the boundary constraints.

\paragraph{Joint deficiency control.}
Under the triangle-tree hypothesis of Section~\ref{sec:joint}, a sufficient condition for~\eqref{eq:target} is a feasible four-coloring for which
\[
\sum_{i=1}^4\Def(R_{I_i})\le
5+\frac{s}{2}+\frac{3D_B}{4}-\frac{5Q}{4}.
\]
For $Q=0$, this is the condition following~\eqref{eq:four-def}. The individual estimate $\Def(R_I)\le P_{R_I}(V(R_I))/3$ does not establish this joint bound: the total degree deficit can grow with the number of claw centers. Moreover, if the displayed right-hand side is negative, no coloring can meet this particular sufficient condition, and a different selection or a stronger use of~\eqref{eq:joint} is needed.

A possible next step is to study the odd-component barriers of the residual roots together with their attachment vertices in $G$. Lemma~\ref{lem:residual-third} rules out odd regions with attachment budget one. Higher-budget regions must be related to selectable claw centers rather than counted as unrelated losses. In the presence of branches, the terminal antichain, the connected exchanges of Lemma~\ref{lem:exchange}, and the one-unit criticality of Proposition~\ref{prop:critical} constrain such barriers simultaneously. Corollary~\ref{cor:neutral} gives a particularly concrete condition: restoring the triangle nodes of any one selected terminal leaves the matching number unchanged and forces the slack in~\eqref{eq:neutral-slack}.

\paragraph{Degree-three components beyond triangle trees.}
Theorem~\ref{thm:joint} assumes that every component of $G[V_3]$ is a triangle tree. A general minimum counterexample has not been shown to satisfy this assumption. Extending the approach therefore also requires boundary profiles for other degree-three components, or structural reductions excluding them. The incidence-forest restriction in Proposition~\ref{prop:incidence} provides a necessary condition for such an extension.

These two tasks delimit the present reductions. The sharp subclass theorem and the residual compatibility theorem are unconditional, but neither the joint deficiency estimate nor the reduction to triangle-tree components is established here.

\clearpage
\appendix
\section{Illustrative triangulations}\label{app:illustrations}
\begin{figure}[!ht]
\centering
\begin{minipage}[t]{.48\textwidth}
\centering
\begin{tikzpicture}[line cap=round,line join=round,
  vertex/.style={circle,draw=black,fill=white,inner sep=0pt,minimum size=13pt,font=\scriptsize},
  low/.style={vertex,double,double distance=.8pt}]
\coordinate (p0) at (0.00000,3.00000);
\coordinate (p1) at (-2.80000,-1.60000);
\coordinate (p2) at (-0.67586,0.40055);
\coordinate (p3) at (0.00000,0.98670);
\coordinate (p4) at (0.67586,0.40055);
\coordinate (p5) at (2.80000,-1.60000);
\coordinate (p6) at (0.00000,-1.06482);
\coordinate (p7) at (-0.57931,-0.52964);
\coordinate (p8) at (0.00000,0.14571);
\coordinate (p9) at (0.57931,-0.52964);
\draw[line width=.45pt] (p0)--(p1);
\draw[line width=.45pt] (p0)--(p2);
\draw[line width=.45pt] (p0)--(p3);
\draw[line width=.45pt] (p0)--(p4);
\draw[line width=.45pt] (p0)--(p5);
\draw[line width=.45pt] (p1)--(p2);
\draw[line width=.45pt] (p1)--(p5);
\draw[line width=.45pt] (p1)--(p6);
\draw[line width=.45pt] (p1)--(p7);
\draw[line width=.45pt] (p2)--(p3);
\draw[line width=.45pt] (p2)--(p7);
\draw[line width=.45pt] (p2)--(p8);
\draw[line width=.45pt] (p3)--(p4);
\draw[line width=.45pt] (p3)--(p8);
\draw[line width=.45pt] (p4)--(p5);
\draw[line width=.45pt] (p4)--(p8);
\draw[line width=.45pt] (p4)--(p9);
\draw[line width=.45pt] (p5)--(p6);
\draw[line width=.45pt] (p5)--(p9);
\draw[line width=.45pt] (p6)--(p7);
\draw[line width=.45pt] (p6)--(p9);
\draw[line width=.45pt] (p7)--(p8);
\draw[line width=.45pt] (p7)--(p9);
\draw[line width=.45pt] (p8)--(p9);
\node[vertex] at (p0) {0};
\node[vertex] at (p1) {1};
\node[vertex] at (p2) {2};
\node[low] at (p3) {3};
\node[vertex] at (p4) {4};
\node[vertex] at (p5) {5};
\node[low] at (p6) {6};
\node[vertex] at (p7) {7};
\node[vertex] at (p8) {8};
\node[vertex] at (p9) {9};
\end{tikzpicture}

\small (a) $4^2 5^8$: vertices $3$ and $6$ are nonadjacent.
\end{minipage}\hfill
\begin{minipage}[t]{.48\textwidth}
\centering
\begin{tikzpicture}[line cap=round,line join=round,
  vertex/.style={circle,draw=black,fill=white,inner sep=0pt,minimum size=13pt,font=\scriptsize},
  low/.style={vertex,double,double distance=.8pt}]
\coordinate (p0) at (-0.29903,0.89801);
\coordinate (p1) at (-0.00000,0.17704);
\coordinate (p2) at (0.29903,0.89801);
\coordinate (p3) at (0.00000,3.00000);
\coordinate (p4) at (-0.97864,0.54260);
\coordinate (p5) at (-0.81553,-0.12759);
\coordinate (p6) at (-0.00000,-0.65563);
\coordinate (p7) at (0.81553,-0.12759);
\coordinate (p8) at (0.97864,0.54260);
\coordinate (p9) at (2.80000,-1.60000);
\coordinate (p10) at (-2.80000,-1.60000);
\draw[line width=.45pt] (p0)--(p1);
\draw[line width=.45pt] (p0)--(p2);
\draw[line width=.45pt] (p0)--(p3);
\draw[line width=.45pt] (p0)--(p4);
\draw[line width=.45pt] (p0)--(p5);
\draw[line width=.45pt] (p1)--(p2);
\draw[line width=.45pt] (p1)--(p5);
\draw[line width=.45pt] (p1)--(p6);
\draw[line width=.45pt] (p1)--(p7);
\draw[line width=.45pt] (p2)--(p3);
\draw[line width=.45pt] (p2)--(p7);
\draw[line width=.45pt] (p2)--(p8);
\draw[line width=.45pt] (p3)--(p4);
\draw[line width=.45pt] (p3)--(p8);
\draw[line width=.45pt] (p3)--(p9);
\draw[line width=.45pt] (p3)--(p10);
\draw[line width=.45pt] (p4)--(p5);
\draw[line width=.45pt] (p4)--(p10);
\draw[line width=.45pt] (p5)--(p6);
\draw[line width=.45pt] (p5)--(p10);
\draw[line width=.45pt] (p6)--(p7);
\draw[line width=.45pt] (p6)--(p9);
\draw[line width=.45pt] (p6)--(p10);
\draw[line width=.45pt] (p7)--(p8);
\draw[line width=.45pt] (p7)--(p9);
\draw[line width=.45pt] (p8)--(p9);
\draw[line width=.45pt] (p9)--(p10);
\node[vertex] at (p0) {0};
\node[vertex] at (p1) {1};
\node[vertex] at (p2) {2};
\node[vertex] at (p3) {3};
\node[low] at (p4) {4};
\node[vertex] at (p5) {5};
\node[vertex] at (p6) {6};
\node[vertex] at (p7) {7};
\node[low] at (p8) {8};
\node[vertex] at (p9) {9};
\node[vertex] at (p10) {10};
\end{tikzpicture}

\small (b) $4^2 5^8 6^1$: vertices $4$ and $8$ are nonadjacent.
\end{minipage}
\caption{Two simple plane triangulations with nonadjacent degree-four
vertices, marked by double circles. Each drawing includes every vertex
and edge; the degree sequences can be read directly from the drawings.
These are examples with \emph{nonadjacent} degree-four vertices, not
counterexamples to the exclusions in Lemma~\ref{lem:finite}.}
\label{fig:triangulation-examples}
\end{figure}
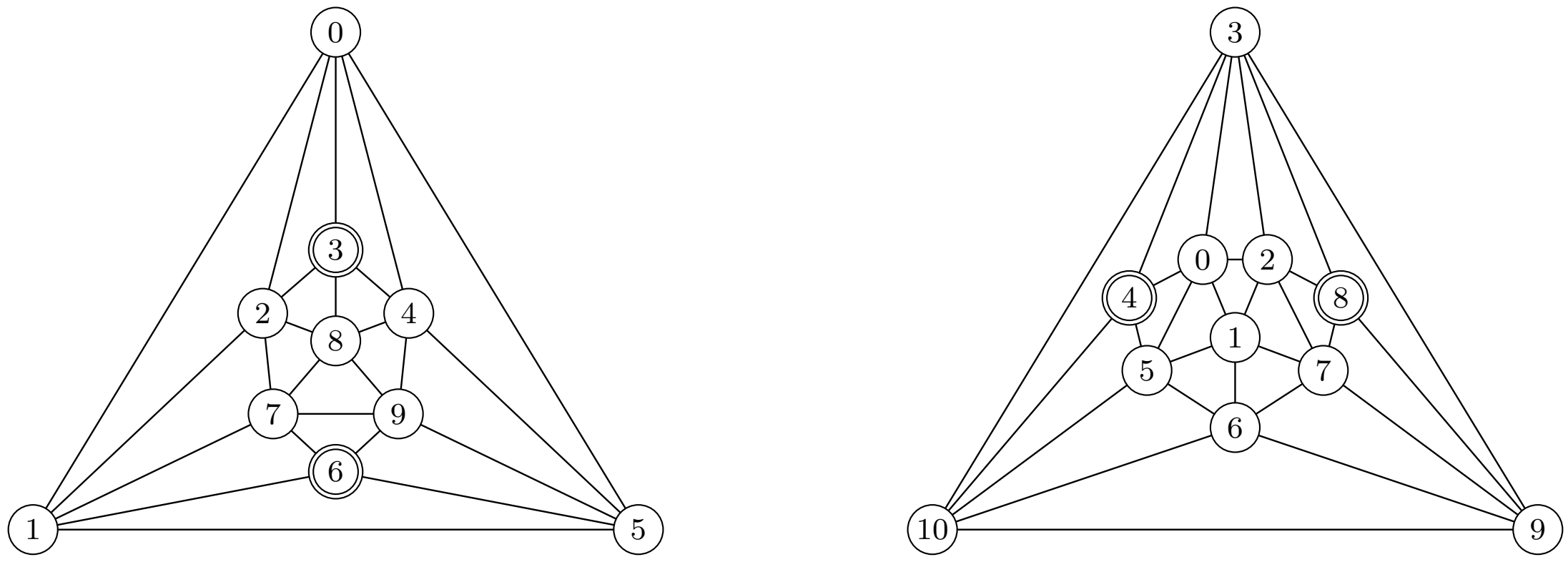

The adjacency condition in Lemma~\ref{lem:finite} is essential.
Figure~\ref{fig:triangulation-examples} exhibits simple plane triangulations
with degree sequences $4^2 5^8$ and $4^2 5^8 6^1$, respectively, whose
degree-four vertices are \emph{not} adjacent. In the first drawing these
vertices are $3$ and $6$; in the second they are $4$ and $8$.
All vertices and edges are shown, so the degrees and nonadjacencies can be
checked directly. These illustrations make no completeness or uniqueness
claim and are not used in the analytic exclusion proof.

\section*{Acknowledgements}
This project is  supported by the National Natural Science Foundation of China (No. 12271246). 
\section*{AI Declaration}

We used ChatGPT to assist with the exploration of proof strategies and the preparation of this manuscript.  All mathematical arguments and proofs in the final manuscript were  verified by the author.

\end{document}